\documentclass[%
 reprint,
showpacs,
amsmath,amssymb,
aps,
pre,
]{revtex4-1}

\usepackage{graphicx}
\usepackage{dcolumn}
\usepackage{bm}
\usepackage{caption}
\usepackage{subfig}
\usepackage{ dsfont }

\graphicspath{{figures/}}
 \newcommand{\ee}{{\rm e}\hspace{1pt}}
\newcommand{\norm}[1]{\Vert #1 \Vert}
\newcommand{\abs}[1]{\left| #1 \right|}
\newcommand{\dd}{{\rm d}}

\newcommand{\trans}{{\sf T}}

\begin{document}

\preprint{APS/123-QED}

\title{Splitting methods for time integration of trajectories in combined electric and magnetic fields}
\thanks{The work was financially supported by the FWF doctoral program `Computational Interdisciplinary Modelling' W1227.
The first author was additionally supported by the FWF project Y398.}%

\author{Christian Knapp \footnotemark[1]}
\author{Alexander Kendl}%
\affiliation{%
 Institute for Ion Physics and Applied Physics, University of Innsbruck, Technikerstra\ss e 25/3, A-6020 Innsbruck, Austria
}%

\author{Antti Koskela}
\affiliation{%
Department of Mathematics, KTH Royal Institute of Technology, Lindstedtv\"agen 25, 10044 Stockholm, Sweden
}%

  \author{Alexander Ostermann}%
\affiliation{%
 Department of Mathematics, University of Innsbruck, Technikerstra\ss e 13, A-6020 Innsbruck, Austria
}%


\date{\today}

\begin{abstract}
The equations of motion of a single particle subject to an arbitrary electric and
a static magnetic field form a Poisson system. We present a second-order time integration method  
which preserves well the Poisson structure and compare it to commonly used algorithms,
such as the Boris scheme. 
All the methods are represented in a general framework of splitting methods.
We use the so-called $\phi$ functions, which give efficient ways for both analyzing
and implementing the algorithms.
Numerical experiments show an excellent long term stability for the new method considered.

\end{abstract}

\maketitle


\section{Introduction}

The time integration of a particle motion due to electric and magnetic fields is 
a common task in several fields of science, from molecular dynamics to 
plasma physics and astrophysics. Various integration methods have been established 
within different scientific areas. A prominent example is the leap frog 
method, which is also called Boris pusher in the context of particle-in-cell (PIC) codes. 

We consider here time-independent magnetic fields 
which, however, may be nonuniform in space. The electric field on the other 
hand is not subject to any constraints. A classical Coulomb system 
subject to a magnetic field serves as an example. If the particles 
within this system are moving sufficiently slow, the pairwise interaction 
can be described purely electrostatically and only the time invariant external 
magnetic field has to be taken into account.

In this article we review some second-order methods from the literature
under a unifying framework: the methods are interpreted as splitting methods and
the solutions of the subflows are given in terms of the so-called $\phi$ functions.
This offers an economic way of stating, analyzing and implementing the methods.
We also present a new second-order method based on the symmetric splitting of the Hamiltonian 
of the system. Further we discuss 
the structure preserving character of the methods regarding time symmetry, 
volume preservation and preservation of the system's underlying Poisson structure.

The article is organized as follows. In Section~\ref{sec:2} we formulate the
equations of motion of a  particle subject to an electric and
a static magnetic field as a Poisson system, and we state some important properties
of the exact flow of the system. In Section~\ref{sec:3} we shortly discuss the corresponding
properties we desire from the discrete flow given by the numerical integrator.
In Section~\ref{sec:splitting} we give formulas for the maps corresponding
to the subflows of different terms of the Poisson system and shortly discuss
their structural properties. In Section~\ref{sec:methods} we present four time integrators, 
which are constructed using the formulas of Section~\ref{sec:splitting}
and certain $\phi$ functions. In Section~\ref{sec:num_results} we give numerical results
in which the integrators are compared on four different cases.
The needed formulas for the matrix functions are given in the appendix (Section~\ref{sec:appendix}).

\section{Equations of motion and their description as a Poisson system} \label{sec:2}

The equations of motion for a particle with mass $m$ and charge $c$ subject 
to an electromagnetic field are given by the Lorentz force
\begin{equation} \label{eq:Lorentz_force}
  m \ddot{q} = c(e + \dot{q}\times b),
\end{equation}
where the electric field $e$ and the magnetic field $b$ are 
determined by a scalar potential $\Phi$ and a vector potential $A$, respectively,
as
\begin{equation*} 
\begin{aligned}
  b &= \nabla \times A, \\
  e &= - \nabla \Phi - \tfrac{\partial}{\partial t}A.
\end{aligned}
\end{equation*}
In the Lagrangian formalism the Lorentz force \eqref{eq:Lorentz_force} 
can be derived from the Lagrangian
\begin{equation} \label{eq:Lagrangian_vecpot}
  L = \tfrac{1}{2} m \norm{\dot{q}}^2 + c\, \dot{q} \cdot A - c\, \Phi,
\end{equation}
which leads to the usual Hamiltonian formulation
\begin{equation*}
\begin{aligned}
  H &= \tfrac{1}{2m} \norm{p - c\, A}^2 + c\, \Phi,  \\
  p &= m \dot{q} + c\,A.
\end{aligned}
\end{equation*}
However, we use a different formulation of the equations of motion, which
exhibits both the Hamiltonian structure of the system and the splitting of the right-hand side
into efficiently computable subflows. We use $p = m \dot{q}$ for the momentum
and state the equations of motion as
\begin{equation} \label{ODE}
\begin{aligned}
   \frac{\dd}{\dd t} q &= \frac{p}{m}, \\
  \frac{\dd}{\dd t} p &= F(q) + \Omega(q) p,
\end{aligned}
\end{equation}
where $F(q)=c\,e(q)=-c\,\nabla \Phi(q)$ is the electric force and
\begin{equation} \label{def:Omega}
\Omega(q) = \begin{bmatrix} 0 & \omega_3(q) & -\omega_2(q) \\ -\omega_3(q) 
& 0 & \omega_1(q) \\ \omega_2(q) & -\omega_1(q) & 0 \end{bmatrix}
\end{equation}
corresponds to the cross product of the magnetic force as
$$
\Omega(q) p = p \times \omega(q) = c\, \dot{q} \times b(q).
$$ 
We use the matrix formulation of the cross product as it will simplify the analysis and derivation
of the time integration methods.

The system \eqref{ODE} can be described as a Poisson system~\cite[Ch.\;7]{HairerLubichWanner}. Recall that a system of ordinary differential equations (ODEs) in 
$\mathbb{R}^n$ with a given initial value $y_0$ at $0$,
\begin{equation} \label{def:Poisson_system}
\frac{\dd}{\dd t} y = B(y) \nabla_y H(y), \quad y(0) = y_0,
\end{equation}
is called a Poisson system, if $B(y)$ is a skew-symmetric matrix which satisfies the identity
\begin{equation} \label{eq:Jacobi_indentity}
\begin{aligned}
\sum\limits_{\ell=1}^n  &\Big( \frac{\partial B_{ij}(y)}{ \partial y_\ell} B_{\ell k}(y) +  
\frac{\partial B_{jk}(y)}{\partial y_\ell} B_{\ell i}(y) \\ 
  & + \frac{\partial B_{ki}(y)}{ \partial y_\ell} B_{\ell j} (y) \Big) = 0 \quad \textrm{ for } i,j,k=1,...,n.
\end{aligned}
\end{equation}
In this case, $B(y)$ is also called a Poisson tensor. The exact flow $\varphi_t: \, \, y(0) \mapsto y(t)$ 
of any Poisson system is a Poisson map 
(also called noncanonically symplectic map), i.e., the Jacobian matrix of the flow satisfies
\begin{equation} \label{eq:Poisson_condition}
 \left( \frac{\partial \varphi_t(y)}{\partial y} \right)  B(y) 
 \left( \frac{\partial \varphi_t(y)}{\partial y}\right)^\trans = B \big(\varphi_t(y) \big). 
\end{equation}
Moreover, it is time symmetric, i.e. $\varphi_{-h} = \varphi_h^{-1}$. 

The ODE \eqref{ODE} can be written formally in the form \eqref{def:Poisson_system} with the Hamiltonian
\begin{equation} \label{eq:Hamiltonian}
H(p,q) = \frac{1}{2m} \norm{p}^2 + c\, \Phi(q)
\end{equation}
and the skew-symmetric structure matrix
\begin{equation} \label{def:system_matrix}
B(y) = \begin{bmatrix} 0 & I \\ -I & m\Omega(q) \end{bmatrix}
\end{equation}
by setting $y = \begin{bmatrix} q \\ p \end{bmatrix}$. We will denote from now on the conjugate variables
as $y$ or $\begin{bmatrix} q \\ p \end{bmatrix}$, whichever is more convenient.

For the system matrix \eqref{def:system_matrix} we have the following.
Let $b:\mathbb{R}^3 \rightarrow \mathbb{R}^3, \quad q \mapsto b(q) = \begin{bmatrix} b_1(q) & b_2(q) & b_3(q) \end{bmatrix}^\trans$ 
be a differentiable function, and let the $3 \times 3$ matrix \eqref{def:Omega}
be defined 
by $\omega(q) = c\, m^{-1}b(q)$. 
Then a direct calculation using the relation
\eqref{eq:Jacobi_indentity} shows that $B(y)$ as in \eqref{def:system_matrix} 
defines a Poisson tensor if and only if $\,\, \nabla \cdot b(q) = 0$.
This is satisfied as Maxwell's laws are assumed to hold for the magnetic field
(e.g., we do not consider magnetic monopoles).
Thus the ODE \eqref{ODE} can be written as a Poisson system \eqref{def:Poisson_system}.


The flow of a Hamiltonian system is volume preserving, i.e., it holds that
\begin{equation} \label{eq:det_ham}
\det \,  \frac{\partial \varphi_t(y)}{\partial y} = 1.
\end{equation}
The same is true for the exact flow of the system \eqref{def:Poisson_system} 
with the structure matrix \eqref{def:system_matrix} as can be seen as follows.
From the relation \eqref{eq:Poisson_condition} it follows that
$$
\abs{\det \,  \frac{\partial \varphi_t(y)}{\partial y}  } = 1,
$$
as $\det \, B(y) = 1$. Since $\det \,   \frac{\partial \varphi_t(y)}{\partial y} $
is a continuous function of $t$, and since $\, \frac{\partial \varphi_0(y)}{\partial y} = I$,
\eqref{eq:det_ham} holds also for the exact flow of the Poisson system.

As we consider only the motion of a single particle, the electrostatic force $F(q)$ will be time-independent.
In case of many particles and pairwise forces this is not the case. However, the Poisson system
\eqref{ODE} generalizes from dimension $6$ to $6N$ in a trivial way: $\Omega(q)$ is replaced
by a $3N \times 3N$ block diagonal matrix consisting of $3 \times 3$ blocks.

\section{Structure preserving integrators} \label{sec:3}

It is desirable that the one-step map $\varphi_h: \,\,  y_0 \mapsto y_1$ given by the numerical flow
satisfies the same structural properties as the exact flow of the ODE.
The preservation of time symmetry or the Poisson structure 
implies good long time behavior for the numerical solution (see the backward
error analysis of \cite{HairerLubichWanner}). 
The properties we want from the integrators are:
\begin{enumerate}
 \item Time symmetry: $\varphi_{-h} = (\varphi_h)^{-1}$.
 \item The integrator gives a Poisson map, i.e., it holds that
     \begin{equation} \label{eq:Poisson_condition2}
 \left( \frac{\partial y_1}{\partial y_0} \right)  B(y_0) \left( \frac{\partial y_1}{\partial y_0} \right)^\trans = B(y_1).       
      \end{equation}
  \item The preservation of volume:
      $$
      \det \,  \frac{\partial y_1}{\partial y_0} = 1.
      $$	
\end{enumerate}
Whether \eqref{eq:Poisson_condition2} is satisfied depends both on the integrator and the Poisson system
considered. Moreover, the integrator has to respect the Casimirs (the first integrals of the Poisson system), 
which in general is difficult to verify for a nonlinear system matrix $B(y)$.
      
However, in some special cases one can show that these properties hold. That is the case
for example if the integrator is a composition method resulting from a Hamiltonian splitting.
This will be discussed in Section~\ref{sec:splitting}.

We will study analytically in Section~\ref{sec:methods} whether the considered integrators satisfy these properties.

\section{Splitting methods} \label{sec:splitting}
To have a common description for the various time integrators we state them
in the framework of splitting methods and use the $\phi$ functions for representing the subflows
containing the term $\Omega(q)p$.
To this end, we write the Poisson system \eqref{ODE} as
\begin{equation*} 
 \frac{\dd }{\dd t} y  = 
\mathcal{T}(y) + \mathcal{E}(y) + \mathcal{B}(y),
\end{equation*}
where
\begin{equation*} 
\begin{aligned}
\mathcal{T} \left( \begin{bmatrix} q \\ p \end{bmatrix} \right) &= \begin{bmatrix} \frac{1}{m}p \\  0 \end{bmatrix}, \\
\mathcal{E}\left( \begin{bmatrix} q \\ p \end{bmatrix} \right) &= \begin{bmatrix} 0 \\  F(q) \end{bmatrix}, \\
\mathcal{B}\left( \begin{bmatrix} q \\ p \end{bmatrix} \right) &= \begin{bmatrix} 0 \\ \Omega ( q ) p \end{bmatrix}. \\
\end{aligned}
\end{equation*}
We denote the exact flows of the subsystems corresponding to $\mathcal{T}$, $\mathcal{B}$, 
$\mathcal{T}+\mathcal{B}$, etc. by
$\varphi_t^\mathcal{T}$, $\varphi_t^\mathcal{B}$, 
$\varphi_t^{\mathcal{T} + \mathcal{B}}$, etc. This means that $\varphi_t^\mathcal{T}(y_0)$
provides the exact solution of the system $w'(t) = \mathcal{T} (w(t)), \, w(0) = y_0$ at time $t$,
and so on.
We further denote by $\varphi_h^{\mathcal{B}(y^*)}$ the exact flow of the ODE
with right-hand side $\mathcal{B}(y^*)$, where 
$y^* = \begin{bmatrix} q^* \\ p \end{bmatrix}$.
This means that the magnetic field is fixed at $q^*$, i.e. $\Omega(q) \equiv \Omega(q^*)$.
One easily verifies the following solutions for the exact flows:
\begin{equation} \label{eq:flows1}
\begin{aligned}
 \varphi_h^\mathcal{T}\left( \begin{bmatrix} q_0 \\ p_0 \end{bmatrix} \right) &= 
\begin{bmatrix} q_0 + \tfrac{h}{m}p_0 \\ p_0 \end{bmatrix},  \\ 
    \varphi_h^\mathcal{E}\left( \begin{bmatrix} q_0 \\ p_0 \end{bmatrix} \right) &= 
    \begin{bmatrix} q_0  \\ p_0 + hF(q_0) \end{bmatrix},  \\
    \varphi_h^{\mathcal{B}(y^*)}\left( \begin{bmatrix} q_0 \\ p_0 \end{bmatrix} \right) &= 
    \begin{bmatrix} q_0 \\ \exp \big(h \Omega(q^*) \big) p_0 \end{bmatrix}.  
\end{aligned}
\end{equation}
We note that the so-called Rodrigues formula allows an efficient implementation 
of the matrix exponential $\exp(h \Omega)$ and the $\phi$ functions
(see~\eqref{eq:Rodrigues} and~\eqref{eq:Rodrigues_phi} of the Appendix).

We will also need the solutions of the subsystems
corresponding to the fields $\mathcal{T} + \mathcal{B}(y^*)$ and $\mathcal{E} + \mathcal{B}(q^*)$.
The flow corresponding to $\mathcal{T} + \mathcal{B}(y^*)$ is given by the exact solution of
the ODE
$$
\frac{\dd }{\dd t}\begin{bmatrix} q \\ p \end{bmatrix} =
\begin{bmatrix} 0 & \tfrac{1}{m}I \\ 0 & \Omega(q^*) \end{bmatrix} \begin{bmatrix} q \\ p \end{bmatrix}, \quad 
$$
which can be expressed by the matrix exponential as
$$
\begin{bmatrix} q(h) \\ p(h) \end{bmatrix} = 
\exp \left( \begin{bmatrix} 0 & \tfrac{1}{m}I \\ 0 & \Omega(q^*) \end{bmatrix} \right) \begin{bmatrix} q_0 \\ p_0 \end{bmatrix}.
$$
Using \eqref{eq:phi_augmented} we see that
\begin{equation} \label{eq:flow_T_plus_B}
\varphi_h^{\mathcal{T} + \mathcal{B}(y^*)} \left( \begin{bmatrix} q_0 \\ p_0 \end{bmatrix} \right)
= \begin{bmatrix} q_0 + \frac{h}{m} \phi_1 \big(h\Omega(q^*)\big) p_0 \\ \exp(h \Omega(q^*)) p_0 \end{bmatrix},
\end{equation}
where $\phi_1(z) = (\exp(z) - 1)/z$. 

On the other hand, the flow corresponding to the field $\mathcal{E} + \mathcal{B}$ is 
given by the exact solution of the ODE
$$
\frac{\dd }{\dd t}\begin{bmatrix} q \\ p \end{bmatrix} =
\begin{bmatrix} 0 \\    \Omega(q)p  +  F(q) \end{bmatrix}.
$$
This solution is given by the variation-of-constants formula \eqref{eq:voc} as
\begin{equation} \label{eq:flow_E_plus_B}
\varphi_h^{\mathcal{E} + \mathcal{B}}\left( \begin{bmatrix} q_0 \\ p_0 \end{bmatrix} \right)
= \begin{bmatrix} q_0 \\ \exp \big(h \Omega(q_0) \big) p_0 + h \phi_1 \big(h \Omega(q_0) \big) F(q_0) \end{bmatrix}.
\end{equation}
We also state the solution for the flow corresponding to $\mathcal{T} + \mathcal{E} + \mathcal{B}$
in the particular case when both $F$ and $\Omega$ are constant. 
This is presented also in \cite{Chin} using another formulation.
The ODE corresponding to the flow is now
$$
\frac{\dd }{\dd t}\begin{bmatrix} q \\ p \end{bmatrix} =
\begin{bmatrix} 0 & \tfrac{1}{m}I \\ 0 & \Omega \end{bmatrix} \begin{bmatrix} q \\ p \end{bmatrix}
+ \begin{bmatrix} 0 \\ F \end{bmatrix},
$$
and its solution is given again by the variation-of-constants formula 
and by formula \eqref{eq:phi_augmented} as
\begin{equation*} \label{eq:flow_last}
  \begin{aligned}
  \varphi_h\left( \begin{bmatrix} q_0 \\ p_0 \end{bmatrix} \right) & = 
  \exp \left( h \begin{bmatrix} 0 & \tfrac{1}{m}I \\ 0 & \Omega \end{bmatrix} \right) 
  \begin{bmatrix} q_0 \\ p_0 \end{bmatrix} \\
  & \qquad + h \phi_1  \left( h \begin{bmatrix} 0 & \tfrac{1}{m}I \\ 0 & \Omega 
 \end{bmatrix} \right) \begin{bmatrix} 0 \\ F \end{bmatrix} \\
  &= \begin{bmatrix} q_0 + \frac{h}{m} \phi_1(h \Omega) p_0 + \frac{h^2}{m} \phi_2(h\Omega)F \\
	\exp(h \Omega) p_0 + h \phi_1(h \Omega) F
      \end{bmatrix},
      \end{aligned}
\end{equation*}
where $\phi_2(z) = (\phi_1(z) - 1)/z$. 

We note that all of these substeps are time symmetric, i.e., $\varphi_{-h}^{\mathcal{T}} = (\varphi_{h}^{\mathcal{T}})^{-1}$,
and so on. 

We also note that the subflows \eqref{eq:flows1} and \eqref{eq:flow_E_plus_B} are volume preserving.
Both the maps $\varphi_h^{\mathcal{E} + \mathcal{B}} : \,\, y_0 \mapsto y_1$
and $\varphi_h^{\mathcal{B}(y^*)} : \,\, y_0 \mapsto y_1$ have the Jacobian
\begin{equation*} 
\frac{\partial y_1}{\partial y_0} = 
                                    \begin{bmatrix} I & 0  \\
                                    \frac{\partial p_1}{\partial q_0} & \exp \big(h \Omega(\widetilde q) \big) \end{bmatrix},
\end{equation*}
where $\widetilde q = q_0$ for $\varphi_h^{\mathcal{E} + \mathcal{B}}$ and 
$\widetilde q = q^*$ for $\varphi_h^{\mathcal{B}(y^*)}$.
Its determinant is given by
$$
\det \, \frac{\partial y_1}{\partial y_0}  = \det \, \exp \big( h \Omega(\widetilde q) \big).
$$
As the Cayley transform of a skew-symmetric matrix is a unitary matrix, 
and as $\det \, \exp( h \Omega(\widetilde q) )$ is a continuous function of $h$ and $\exp(0) = I$, we find for  
the maps $\varphi_h^{\mathcal{E} + \mathcal{B}}$ and $\varphi_h^{\mathcal{B}(y^*)}$ the relation
$$
\det \, \frac{\partial y_1}{\partial y_0}  = 1.
$$
Therefore both are volume preserving. One easily verifies from \eqref{eq:flows1} that
$\varphi^{\mathcal{T}}$ and $\varphi_h^{\mathcal{E}}$ are also volume preserving.

Splitting the Hamiltonian \eqref{eq:Hamiltonian} gives subsystems corresponding
to the right-hand sides  
\begin{equation} \label{eq:sub_hamiltonian1}
B(y) \, \nabla_y \tfrac{\norm{p}^2}{2m} = \mathcal{T}(y) + \mathcal{B}(y)
\end{equation}
and
\begin{equation} \label{eq:sub_hamiltonian2}
B(y) \,	 \nabla_y \, c \Phi(q) = \mathcal{E}(y).
\end{equation}
Both of these constitute Poisson systems of the form \eqref{def:Poisson_system}
with the system matrix \eqref{def:system_matrix} and with Hamiltonians 
$\tfrac{\norm{p}^2}{2m}$ and $c \Phi(q)$, respectively. Therefore, the subflows
$\varphi_h^{\mathcal{T} + \mathcal{B}}$ and $\varphi_h^{\mathcal{E}}$
are Poisson maps, as well as their compositions 
$\varphi_h^{\mathcal{E}} \circ \varphi_h^{\mathcal{T} + \mathcal{B}}$, etc.

\section{Methods} \label{sec:methods}

We will next state four time integration methods. They are formulated using the $\phi$ functions and 
the solutions of the subflows given in Section~\ref{sec:splitting}.
Three of them can be found in the literature: the Boris--Buneman scheme \cite{BirdsallLangdon},
the method of Chin \cite{Chin} and the method of Spreiter and Walter \cite{SpreiterWalter}.
We also state a new method which is nearly symmetric and nearly preserves the Poisson structure (method C below).

\subsection{Boris--Buneman scheme} \label{subsec:Boris}
The so-called Boris pusher \cite{BirdsallLangdon} is undoubtedly the most common particle trajectory integrator 
used for PIC codes. 
It is a second-order accurate centered difference leap frog scheme, given by
\begin{align} 
  \frac{q(t+\tfrac{h}{2})-q(t-\tfrac{h}{2})}{h} &= \frac{p(t)}{m}, \nonumber \\
  \frac{p(t+h)-p(t)}{h} &= F \big( q(t+\tfrac{h}{2}) \big)   \label{eq:leap_frog} \\
   &+ \Omega  \big( q(t+\tfrac{h}{2}) \big)  \left( \frac{p(t+h)+p(t)}{2} \right) \nonumber.
\end{align}
PIC codes usually define the position with an offset of half a time step $h$ relative to the momentum.
In \eqref{eq:leap_frog} $q$ is updated alternatingly with $p$.

To formulate this scheme as a one-step method, 
we denote 
$ q_0 = \big( q(t+\tfrac{h}{2}) + q(t-\tfrac{h}{2}) \big)/2$,
$q_{1/2} = q(t+\tfrac{h}{2})$, $p_0 = p(t)$ and so on.
Using this notation the linear drift operation corresponding to $\mathcal{T}$ can be split as
\begin{equation} \label{Boris_drift}
\begin{aligned}
 q_{1/2} &= q_0 + \tfrac{h}{2m} p_0, \\
  q_1 &=q_{1/2} + \tfrac{h}{2m} p_1.
\end{aligned}
\end{equation}
The momentum change in between these drift operations is defined by the second equation of \eqref{eq:leap_frog}.
Solving it for $p_1$ gives
$$
p_1 = h(I-\tfrac{h}{2}\hat \Omega)^{-1} F( q_{1/2} ) + R(\tfrac{h}{2}\hat \Omega)p_0,
$$
where $R(z)$ is the Cayley transform 
$$
R(z) = (1+z)(1-z)^{-1},
$$
and $\hat \Omega = \Omega(q_{1/2})$.
We note that both the resolvent $(I-\tfrac{h}{2}\hat \Omega)^{-1}$ and $R(\tfrac{h}{2}\hat \Omega)$ 
can be evaluated efficiently (see \eqref{eq:inverse_I-O} and \eqref{rot_Boris} of the Appendix).
Since $R(z) + 1 = 2(1-z)^{-1}$, we may rewrite the momentum step also as
$$
p_1 = \tfrac{h}{2} \big(R(\tfrac{h}{2}\hat \Omega ) + I \big) F( q_{1/2} ) + R(\tfrac{h}{2}\hat \Omega)p_0.
$$
From this one can easily read off that the mid-step can be split as
\begin{equation} \label{Boris}
\begin{aligned}
  p^+ &= p_0 + \tfrac{h}{2} F ( q_{1/2} ), \\
  p^{++} &= R(\tfrac{h}{2}\hat \Omega) p^+, \\
  p_1 &= p^{++} + \tfrac{h}{2} F ( q_{1/2} ).
\end{aligned}
\end{equation}
The steps \eqref{Boris_drift} and \eqref{Boris} give together the approximation
\begin{equation} \label{eq:Boris_splitted}
\varphi_h \approx \varphi_{\frac{h}{2}}^{ \mathcal{T}} \circ \varphi_{\frac{h}{2}}^{ \mathcal{E}} \circ
\widehat \varphi_h^{\mathcal{B}(y_{1/2})} \circ \varphi_{\frac{h}{2}}^{ \mathcal{E}} \circ \varphi_{\frac{h}{2}}^{ \mathcal{T}},
\end{equation}
where $y_{1/2} =  \varphi_{\frac{h}{2}}^{ \mathcal{E}} \circ \varphi_{\frac{h}{2}}^{ \mathcal{T}}(y_0)$, and
$\widehat \varphi_h^{\mathcal{B}(y_{1/2})}$ is given by
$$
\widehat \varphi_h^{\mathcal{B}(y_{1/2})} \left( \begin{bmatrix} q \\ p \end{bmatrix} \right) = 
    \begin{bmatrix} q \\ R \big( \tfrac{h}{2} \Omega(q_{1/2}) \big) p \end{bmatrix}. 
$$
Since $R(-z) = R(z)^{-1}$, we have 
$\widehat \varphi_{-h}^{\mathcal{B}(y_{1/2})} = \left(\widehat \varphi_h^{\mathcal{B}(y_{1/2})} \right)^{-1}$. 
Therefore, \eqref{eq:Boris_splitted} gives a symmetric method.
We note that \eqref{eq:Boris_splitted} is equivalent to the formulation \eqref{eq:leap_frog} 
with the initial value $\varphi_{h/2}^{\mathcal{T}}(y_0)$.


Since the Cayley transform of a skew-symmetric matrix is a unitary matrix,
the subflow $\widehat \varphi_h^{\mathcal{B}(y_{1/2})}$ is also a volume preserving map
(see the analysis of Section~\ref{sec:splitting}).
And as also the maps $\varphi_{\frac{h}{2}}^{ \mathcal{T}}$ and $\varphi_{\frac{h}{2}}^{ \mathcal{E}}$
are volume preserving, the method \eqref{eq:Boris_splitted} is volume preserving
as a composition of volume preserving maps. The volume preservation property
of the Boris scheme was also shown in \cite{Qinetal}.

%

Replacing the Cayley transform by the exponential function in the mid-step, i.e., considering the splitting
\begin{equation} \label{eq:Boris_splitted2}
\varphi_h \approx \varphi_{\frac{h}{2}}^{ \mathcal{T}} \circ \varphi_{\frac{h}{2}}^{ \mathcal{E}} \circ
\varphi_h^{\mathcal{B}(y_{1/2})} \circ \varphi_{\frac{h}{2}}^{ \mathcal{E}} \circ \varphi_{\frac{h}{2}}^{ \mathcal{T}}
\end{equation}
gives also a symmetric and volume preserving method.  
Both the Cayley transform and the 
matrix exponential rotate the momentum around the vector
$b(q)$ in the mid-step. This can be seen as follows.
Since $b(q)^\trans \Omega(q) p_0 = \tfrac{c}{m} b(q)^\trans \big( p_0 \times b(q)  \big) =  0 $, 
we may deduce from the power series representations that 
\begin{equation} \label{eq:projections}
b(q)^\trans \exp \big(h \Omega(q) \big) p_0 =  b(q)^\trans R \big(\tfrac{h}{2} \Omega(q)\big) p_0 = b(q)^\trans p_0.
\end{equation}
On the other hand, since $\Omega(q)$ is skew-symmetric both $\exp \big(h \Omega \big)$ and $R \big(\tfrac{h}{2} \Omega(q) \big)$ 
are unitary matrices, so that 
\begin{equation} \label{eq:norm_1}
\norm{\exp \big(h \Omega(q) \big) p_0}_2 = \norm{R \big(\tfrac{h}{2} \Omega(q) \big) p_0}_2 = \norm{p_0}_2.
\end{equation}
From \eqref{eq:projections} and \eqref{eq:norm_1} it follows that both give rotations around $b(q)$.

\subsection{Chin's schemes}
Using the Lie operator formalism \cite{HairerLubichWanner}, 
the following splitting methods are derived in \cite{Chin}:
\begin{itemize}
\item Chin-a:
\begin{equation} \label{eq:China}
\varphi_h^a (y_0) = \varphi_{\frac{h}{2}}^{ \mathcal{E} + \mathcal{B}} \circ \varphi_h^{ \mathcal{T} } \circ 
\varphi_{\frac{h}{2}}^{ \mathcal{E} + \mathcal{B} } (y_0),
\end{equation}
\item Chin-b:
\begin{equation} \label{eq:Chinb}
\varphi_h^b (y_0) = \varphi_{\frac{h}{2}}^{\mathcal{T}} \circ \varphi_h^{ \mathcal{E} + \mathcal{B}} \circ 
\varphi_{\frac{h}{2}}^{ \mathcal{T} } (y_0).
\end{equation}
\end{itemize}
Since the substeps are symmetric maps, both \eqref{eq:China} and \eqref{eq:Chinb}
give symmetric methods. Moreover, as the methods \eqref{eq:China} 
and \eqref{eq:Chinb} are compositions of volume preserving
maps (see Section~\ref{sec:splitting}), both of them are volume preserving as well.

We note that the methods of Chin become exactly energy preserving when $F=0$ \cite{Chin}.

\subsection{Symmetric splitting of the Hamiltonian}

As a new method, we propose a symmetric splitting of the Hamiltonian for the Poisson system \eqref{def:Poisson_system}.
This means that the right-hand side of the \eqref{def:Poisson_system} is split into the
flows \eqref{eq:sub_hamiltonian1} and \eqref{eq:sub_hamiltonian2} according to the Hamiltonian \eqref{eq:Hamiltonian}.
Then, Strang splitting is applied to this 
decomposition giving the symmetric one-step method
\begin{equation} \label{eq:symm_split}
y_1 =  \varphi_{\frac{h}{2}}^{\mathcal{E}} \circ \varphi_h^{\mathcal{T} + \mathcal{B}} \circ \varphi_{\frac{h}{2}}^{\mathcal{E}} (y_0).
\end{equation}
For a constant magnetic field $\Omega(q) \equiv \Omega$ this splitting approach can be found in \cite{LeimkuhlerReich}
where it goes by the name of Scovel's method; see also \cite{PatacchiniHutchinson}.


The mid-step of \eqref{eq:symm_split}, i.e., the solution
\begin{equation} \label{eq:mid_step}
v_1= \varphi_h^{\mathcal{T} + \mathcal{B}}(v_0), \quad \textrm{where} \quad v_0 = \varphi_{h/2}^{\mathcal{E}} (y_0),
\end{equation}
can be computed efficiently using the formula
\eqref{eq:flow_T_plus_B} if the magnetic field is uniform. If the 
magnetic field is position dependent, some approximations have to be made.

We infer from \eqref{eq:symm_split} that the last substep $\varphi_{h/2}^{\mathcal{E}}$ changes only $p$.
As the magnetic field depends only on $q$,
we see from \eqref{eq:flow_T_plus_B} that by using a symmetric approximation for \eqref{eq:mid_step},
we get a symmetric approximation for \eqref{eq:symm_split}.

The approaches we consider for the mid-step are the following ones:
\begin{enumerate}
 \item The symmetric method
  \begin{equation} \label{eq:subsplit1}
w_1 = \varphi_{\frac{h}{2}}^{\mathcal{T} + \mathcal{B}(w_1)} \circ \varphi_{\frac{h}{2}}^{ \mathcal{T} + \mathcal{B}(w_0)}(w_0)
\end{equation}
    combined with fixed-point iteration. This approach can be found in \cite{EinkemmerOstermann}.
 \item 
   The symmetric method
   \begin{equation} \label{eq:subsplit2}
w_1 = \varphi_h^{\mathcal{T} + \mathcal{B}((w_1+w_0)/2)}(w_0)
\end{equation}
   combined with fixed-point iteration. 
  \item
   Perform a symmetric high-order composition of either \eqref{eq:subsplit1} or \eqref{eq:subsplit2}.
   For a description of composition schemes we refer to \cite{HairerLubichWanner}.
   In numerical experiments we use the 8th-order scheme from \cite{SuzukiUmeno}.
  
\end{enumerate}
In the methods above the numerical strategy is to 
perform a fixed-point iteration for the implicit mid-step. 
This approach is also taken in \cite{EinkemmerOstermann}.
For the approach \eqref{eq:subsplit1} this means that after the initial value 
$$
w_1^{(0)} = \varphi_{\frac{h}{2}}^{\mathcal{T} + \mathcal{B}(w_0)} (w_0), \quad w_0 = \varphi_{\frac{h}{2}}^{\mathcal{E}} ( y_0 )
$$ 
is computed, the iteration goes on as 
$$
w_1^{(k)} = \varphi_{\frac{h}{2}}^{\mathcal{T} + \mathcal{B} (w_1^{(k-1)} )} (w_1^{(0)}).
$$
Here $\varphi_{\frac{h}{2}}^{\mathcal{T} + \mathcal{B}(w_1^{(k-1)})} (w_1^{(0)})$ can be computed 
explicitly using the formula \eqref{eq:flow_T_plus_B}. 

The implementation for \eqref{eq:subsplit2} goes analogously.
After the initial value 
$$
w_1^{(0)} = \varphi_{\frac{h}{2}}^{\mathcal{E}} ( y_0 )
$$ 
is computed, the iteration goes on as 
$$
w_1^{(k)} = \varphi_{\frac{h}{2}}^{\mathcal{T} + \mathcal{B} ( (w_1^{(k-1)} + w_1^{(0)})/2  )} (w_1^{(0)}).
$$


As shown in \cite{EinkemmerOstermann} the method \eqref{eq:symm_split} with the approach \eqref{eq:subsplit1} will be symmetric
up to an order that equals the number of fixed-point iterations. It is easy to verify that
the same holds true for \eqref{eq:subsplit2}. By the results of \cite{EinkemmerOstermann}
the same holds true for the third approach.

It is easily verified that the high-order composition scheme will also be
a Poisson integrator up to an order of the scheme used for the mid-step 
(i.e., the relation \eqref{eq:Poisson_condition2} will be satisfied up to this order).
Although the number of magnetic field evaluations 
increases considerably, the term $F(q)$ is still evaluated only once per time step,
which may be beneficial in situations where the calculation of 
the electric field is more expensive than that of the magnetic field.

Since the flows of the subsystems resulting from a Hamiltonian splitting are Poisson maps,
so is the composition \eqref{eq:symm_split}; see also~\cite[Ch.\;7]{HairerLubichWanner}.
As the method \eqref{eq:symm_split} is symmetric and a Poisson integrator, it
has a modified Hamiltonian (see the backward error analysis of~\cite{HairerLubichWanner})
and preserves the first integrals for exponentially long times.

The strategies 1-3 above give symmetric (or Poisson) maps up to a high order and
we expect a good preservation of first integrals numerically.

\subsection{The method of Spreiter and Walter}

In \cite{SpreiterWalter} a second-order accurate scheme is proposed which
incorporates a homogeneous magnetic field in $z$-direction, and which overcomes 
the stability condition $h \omega_c \ll 2 \pi$, where $\omega_c$ is the cyclotron frequency.
We repeat here the derivation of this method using the $\phi$ functions which
shortens the derivation considerably.
As this derivation shows, the method can be seen as an exponential version
of the velocity Verlet scheme.
Integrating the Lorentz force law 
$$
\dot{p}(t) = \Omega p(t) + F(q(t))
$$
with respect to $t$ gives an ODE for $q$,
\begin{equation} \label{eq:lorentz_q1}
\dot{q}(t) - \dot{q}(0) = \Omega \big( q(t) - q(0) \big) + \tfrac{1}{m} \int\limits_0^t F(q(s)) \, \dd s.
\end{equation}
Under the assumption that $\tfrac{\dd}{\dd t} F(q(t)) = \mathcal{O}(1)$,
we can approximate the ODE \eqref{eq:lorentz_q1} up to second order as
$$
\dot{q}(t) = \Omega q(t) + \big( \dot{q}(0) - \Omega q(0) \big) + \frac{t}{m} F(q(0) ) + \mathcal{O}(t^2).
$$
By the variation-of-constants formula \eqref{eq:voc} the solution for this ODE at $t=h$ is given by
\begin{align} 
 q(h) & =  \exp(h \Omega) q(0) + h \phi_1( h \Omega) \big( \dot{q}(0) - \Omega q(0) \big) \nonumber \\
      & \quad \quad + \tfrac{h^2}{m} \phi_2(h \Omega) F(q(0)) + \mathcal{O}(h^3) \label{eq:lorentz_q2} \\
      &=  q(0) + \tfrac{h}{m} \phi_1(h\Omega) p(0) + \tfrac{h^2}{m} \phi_2(h \Omega) F(q(0)) + \mathcal{O}(h^3). \nonumber
\end{align}
This gives the time stepping formula for $q$. Integrating the Lorentz force law
and using the variation-of-constants formula \eqref{eq:voc} once more yields
$$
p(h) = \exp(h \Omega) p(0) + \int\limits_0^h \exp \big((h-s) \Omega \big) F(q(s)) \, \dd s.
$$
Using the fact that
$$
F(q(s)) = F(q(0)) + s \dot{F}(q(0)) + \mathcal{O}(s^2)
$$
and
$$
\dot{F}(q(0)) = \tfrac{1}{h}\big( F(q(h)) - F(q(0)) \big) + \mathcal{O}(h^2)
$$
we find that
\begin{equation} \label{eq:lorentz_p1}
\begin{aligned}
p(h) =& \exp(h\Omega) p(0) + h \phi_1(h \Omega) F(q(0)) \\
   &+ h \phi_2(h \Omega) \big( F(q(h)) - F(q(0)) \big)  + \mathcal{O}(h^3).
\end{aligned}
\end{equation}
The approximations \eqref{eq:lorentz_q2} and \eqref{eq:lorentz_p1} give together the one-step
formula of Spreiter and Walter:
\begin{equation*}
 \begin{aligned}
  q_1 &= q_0 + \tfrac{h}{m}\phi_1(h\Omega) p_0 + \tfrac{h^2}{m} \phi_2(h \Omega) F(q_0), \\
  p_1 &= \exp(h\Omega) p_0 + h \phi_1(h \Omega) F(q_0) \\
      & \quad + h \phi_2(h \Omega) \big(F(q_1) - F(q_0) \big).
 \end{aligned}
\end{equation*}
Since $\exp(h \Omega) \rightarrow 0$, $\phi_1(h\Omega) \rightarrow 1$ and 
$\phi_2(h\Omega) \rightarrow \tfrac{1}{2}$ as $\Omega \rightarrow 0$,
we see that in the limiting case $\Omega \rightarrow 0$ the velocity Verlet scheme is obtained.
We note that the method of Spreiter and Walter can be interpreted also as an
\it exponential Taylor method \rm \cite{KoskelaOstermann}.

Despite of its similarity with the velocity Verlet scheme, we have not found any
structure preserving properties for the method of Spreiter and Walter.

\section{Numerical comparison of the integrators} \label{sec:num_results}

In the first experiment, we compare 
the integrators for a particle motion in an inhomogeneous magnetic 
field in the absence of an electric field.
Then, results for experiments including combined electric and magnetic fields are provided.
The motion of a charged particle in a Penning trap \cite{Kretzschmar}
provides a suitable numerical test 
in the case of combined fields.
In this setting the particle trajectory is a periodic orbit, formed by a superposition 
of three harmonic oscillators. In addition to the ideal Penning trap with a uniform 
magnetic field we consider two variations of this setup with static but nonuniform magnetic fields.

\subsection{2d particle motion in an inhomogeneous magnetic field}

In the first experiment we consider an example provided in \cite{Chin}: 
a particle of charge $c=-1$ and mass $m=1$ 
at initial position $q_0=\begin{bmatrix} 1 & 0 & 0 \end{bmatrix}^{\trans}$
with initial velocity $v_0=\begin{bmatrix} 0 & v & 0 \end{bmatrix}^{\trans}$
moves in a magnetic field $b= \begin{bmatrix} 0 & 0 & b_z \end{bmatrix}^{\trans}$,
$b_z=\frac{1}{q_x^2}$. The particle motion is composed of a periodic motion
in the $xy$-plane and a $\nabla B$ drift in $y$-direction with drift velocity $v_d=\frac{v^2}{1+v}$.
The drift reduced trajectory $\tilde{q}(t)=q(t)- \tilde{q}_0(t)$ forms an ellipse centered 
at $\tilde{q}_0(t) = \begin{bmatrix} x_{mid} & v_d \, t & 0 \end{bmatrix}^{\trans}$ with $x_{mid}=\frac{1+v}{1+2v}$.
$\tilde{q}(t)=\begin{bmatrix} \tilde{x}(t) & \tilde{y}(t) & 0 \end{bmatrix}^{\trans}$ is given by
\begin{equation*} 
\begin{aligned}
\tilde{x}(t) &= \frac{v}{1+2v} \, \cos\rho(t), \\
\tilde{y}(t) &= \frac{v}{(1+v)+\sqrt{1+2v}} \, \sin\rho(t).
\end{aligned}
\end{equation*}
One full orbit has the period $P =\frac{2 \pi (1+v)}{(1+2v)^{\frac{3}{2}}}$.

The scalar potential is set $\Phi=0$ and the vector potential $A=\begin{bmatrix} 0 & A_x & 0 \end{bmatrix}^{\trans}$ 
with $A_x=-\frac{1}{q_x}$. The Lagrangian takes the form 
$L = \tfrac{1}{2}   \norm{\dot{q}}^2 + \, \dot{q}_y  \tfrac{1}{q_x}$.
As it is independent of $q_y$, the $y$-component of the generalized momentum, 
\begin{equation} \label{eq:Invariant}
I:=\tfrac{\partial L}{\partial \dot{q}} = m\, \dot{q}_y+\frac{1}{q_x}, 
\end{equation}
gives another integral of motion. 

Since now $F(q)=0$, the flow $\varphi_t^\mathcal{E}$ vanishes and the 
integrators used in the previous example simplify as follows.

\begin{itemize}
\item Boris--Buneman scheme with Cayley transformation:
\begin{equation} \label{eq:Boris--Buneman_noE}
\varphi_h (y_0) = \varphi_{\frac{h}{2}}^{ \mathcal{T}} \circ
\widehat \varphi_h^{\mathcal{B}(y_{1/2})} \circ
\varphi_{\frac{h}{2}}^{ \mathcal{T}}(y_0),
\end{equation}
\item Chin-b (equals Boris--Buneman scheme with matrix exponential):
\begin{equation*}
\varphi_h (y_0) = \varphi_{\frac{h}{2}}^{\mathcal{T}} \circ
\varphi_h^{\mathcal{B}(y_{1/2})} \circ
\varphi_{\frac{h}{2}}^{ \mathcal{T} } (y_0).
\end{equation*}
\end{itemize}


The method of Scovel gave poor results as it is not symmetric when applied to a nonuniform field.
Therefore it was discarded in this experiment.
The symmetric methods \eqref{eq:subsplit1} (denoted as Impl.~Strang) and 
\eqref{eq:subsplit2} (denoted as Impl.~mp.) were used with five (for Impl.~Strang) and six (for Impl.~mp.) 
fixed-point iteration steps.

\begin{figure}[ht!]
   \centering
      \subfloat[Relative error of invariant $I$.]{\label{fig:inhom_mag_field_invariant}
\hspace{-5mm} \includegraphics[width=90mm]{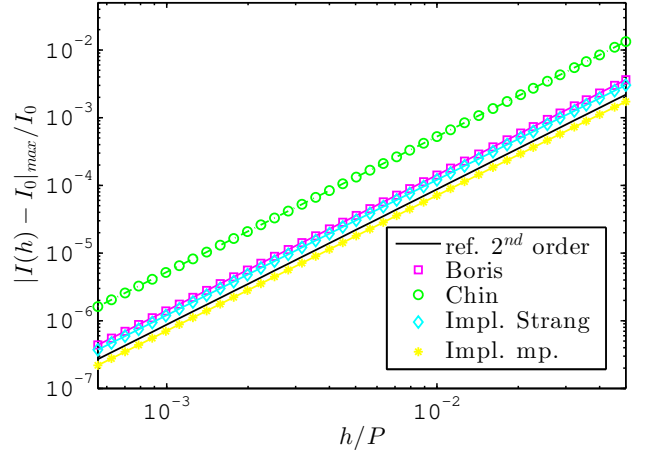}}\\
      \subfloat[Absolute error of position $q$.]{\label{fig:inhom_mag_field_position}
\hspace{-5mm} \includegraphics[width=90mm]{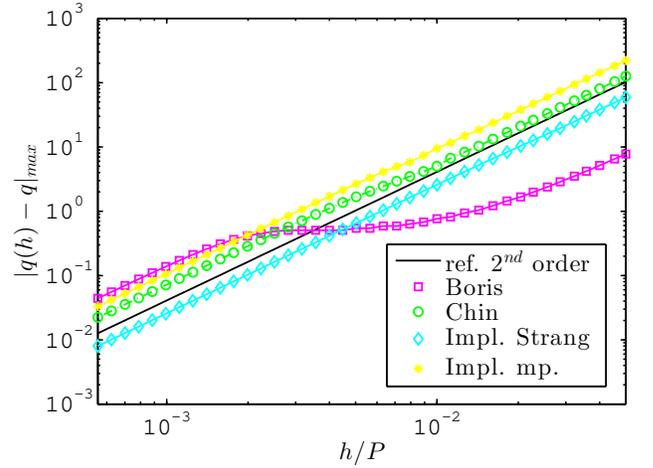}}
   \caption{2d inhomogeneus magnetic field: maximum error (along trajectory) of the invariant 
$I$ and the position $q$ as a function of the time step length $h$ 
(in units of the analytic period length $P$) for a fixed total simulation time of $20000 \, P$.}
\end{figure}

\begin{figure}
   \centering
      \subfloat[Relative error of average drift velocity $\overline{v_d}$.]{\label{fig:inhom_mag_field_v_drift}
\hspace{-5mm} \includegraphics[width=90mm]{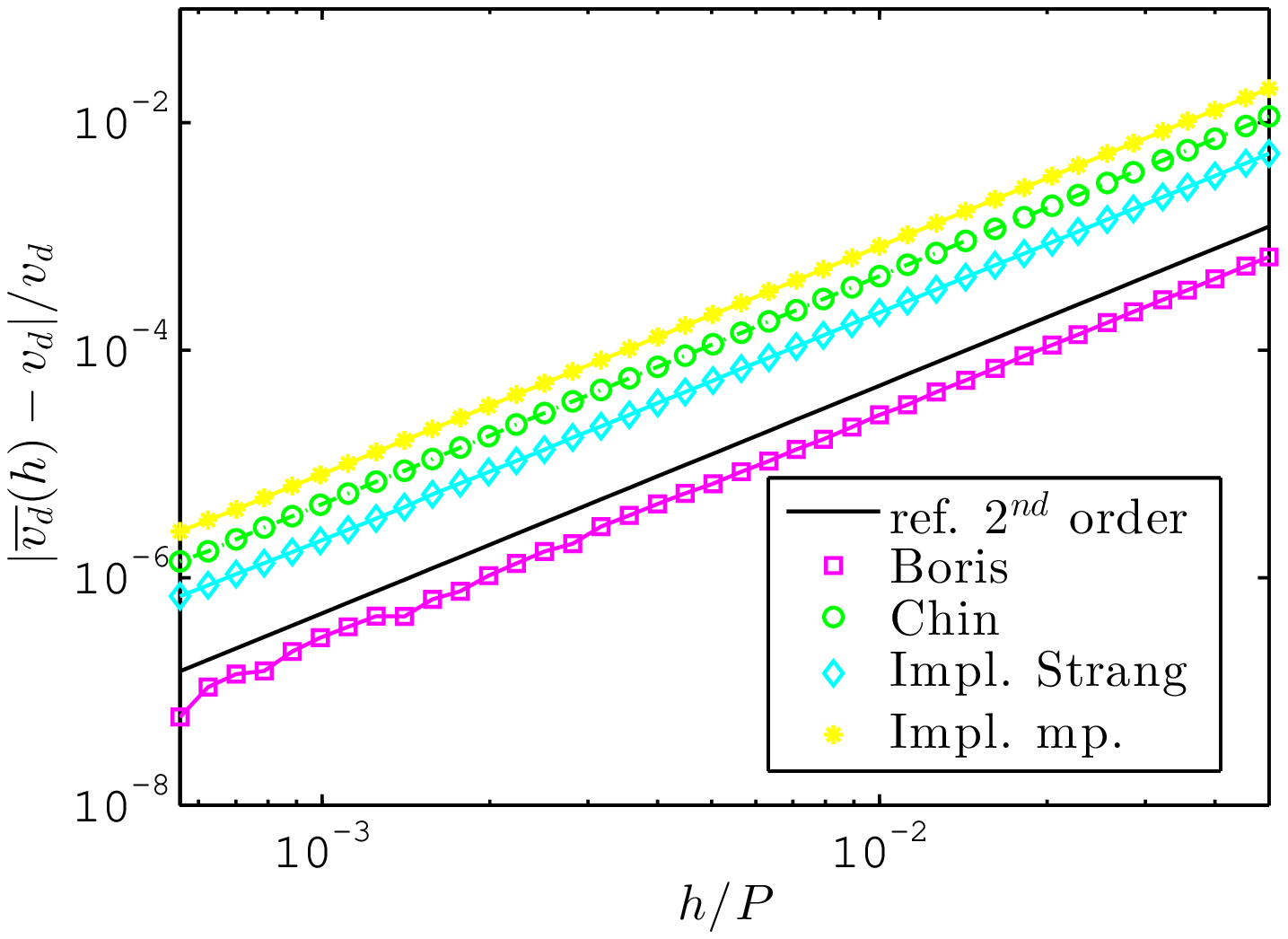}}\\
      \subfloat[Relative error of average period length $\bar{P}$.]{\label{fig:inhom_mag_field_period}
\hspace{-5mm} \includegraphics[width=90mm]{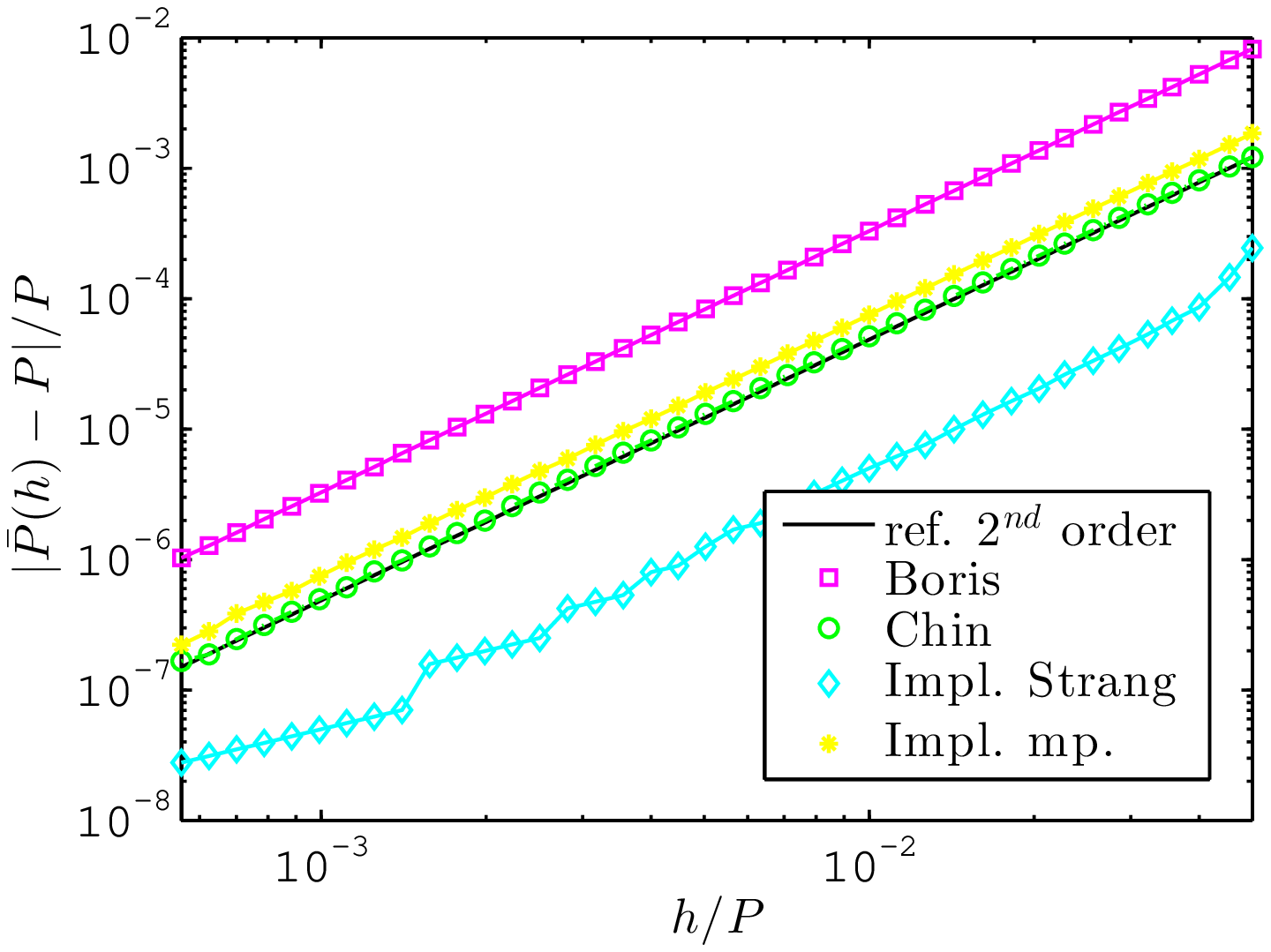}}
   \caption{2d inhomogeneus magnetic field: error of the average drift velocity 
$\overline{v_d}$ and period length $\bar{P}$ as a function of the time step length $h$ 
(in units of the analytic period length $P$) for a fixed total simulation time of $20000 \, P$.}
\end{figure}

We compare these methods by performing a numerical integration 
over 20000 cycles of period $P$ with step sizes $h$ varying from
$0.0005\, P$ to $ 0.05 \, P$. The initial velocity is set to $v=0.5$ in all test cases. 
Figures~\ref{fig:inhom_mag_field_invariant} and~\ref{fig:inhom_mag_field_position} 
show the time step length plotted against the relative error of the invariant $I$ (as defined in \eqref{eq:Invariant})
and the absolute position error using a logarithmic scale.
Both iterative methods provide slightly lower errors of the invariant, while only the 
implicit Strang splitting method has a lower position error than the method of Chin. 
The method of Boris preserves well first-order drifts.
When the step size increases almost about one order of magnitude there is only 
a little impact on the position error for the method of Boris. The position error stays around $\Delta q = 0.5$, 
which is the diameter of the ellipse in $x$--direction.
For small time 
steps, however, the implicit Strang splitting method is the most accurate. 
The error behavior becomes more visible when considering separately the error of the drift 
velocity and the error of the period length as shown in Figures~\ref{fig:inhom_mag_field_v_drift} 
and~\ref{fig:inhom_mag_field_period}. The Boris pusher is the best method regarding 
the conservation of the drift velocity but the worst regarding the period length.

\begin{figure}
   \centering
      \subfloat[Uniform magnetic field.]{ \label{fig:trajectories_a}\includegraphics[width=90mm]{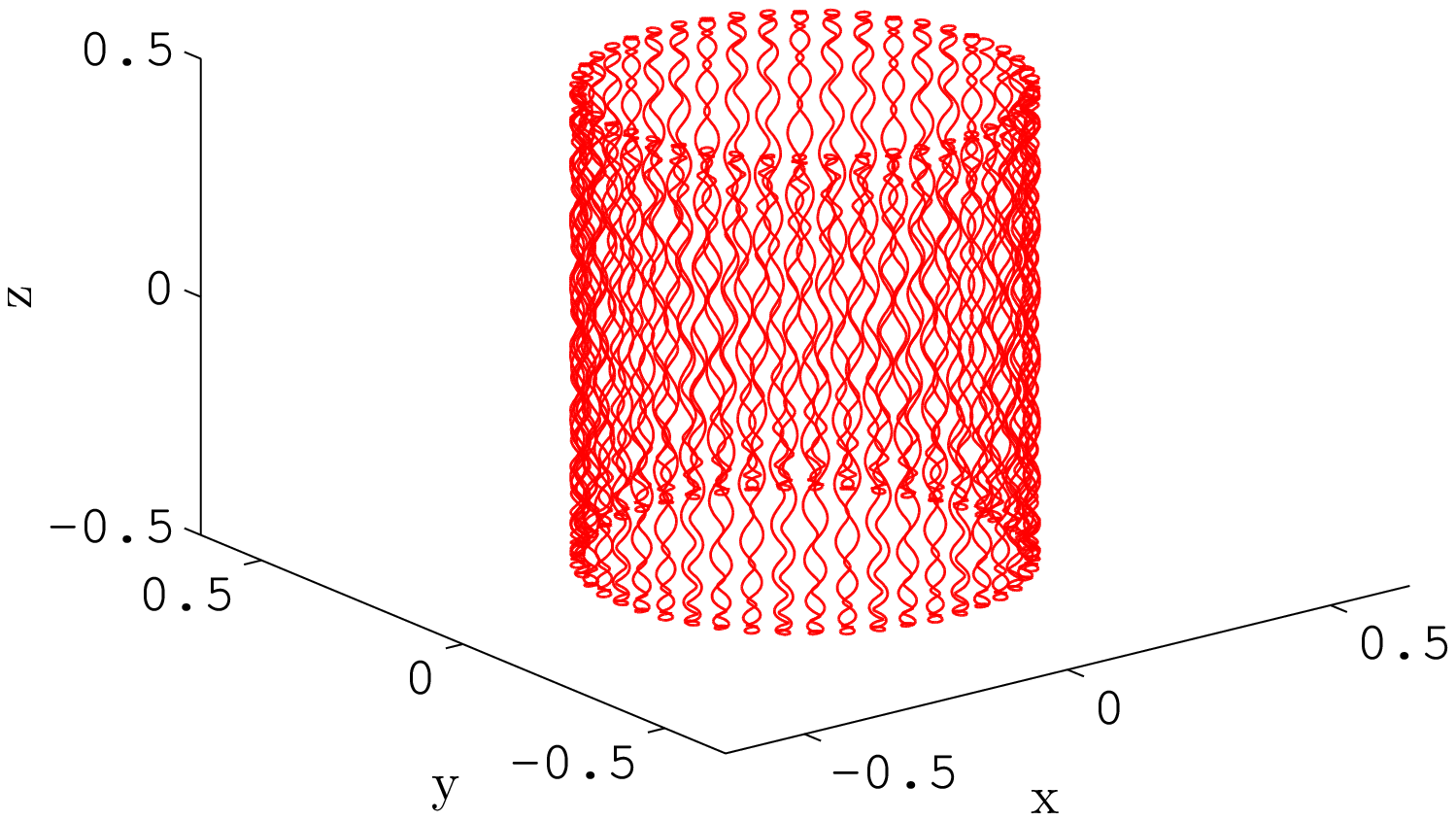}} \\
      \subfloat[Magnetic bottle.]{ \label{fig:trajectories_b} \includegraphics[width=90mm]{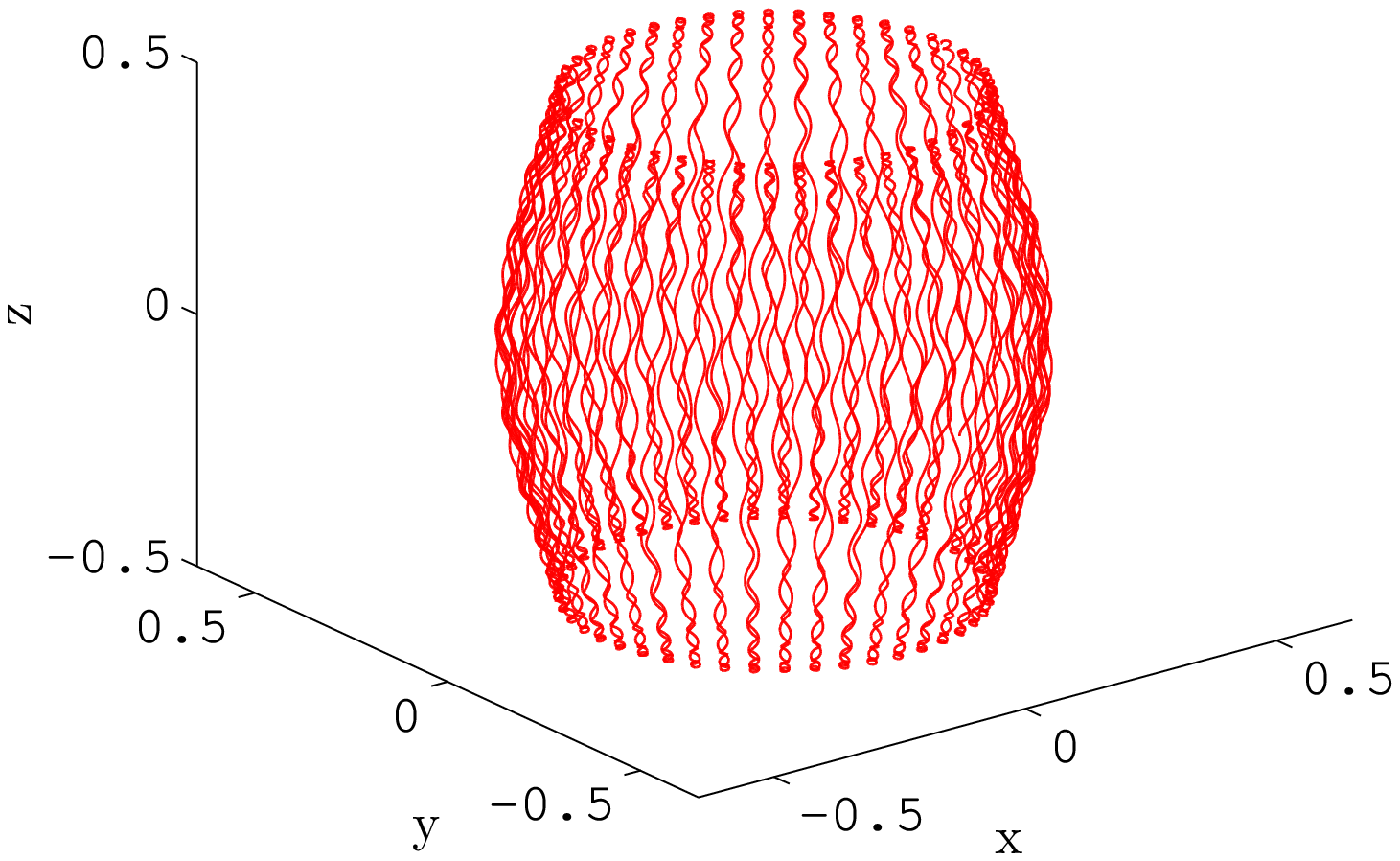}} \\
      \subfloat[Asymmetric magnetic field.]{\label{fig:trajectories_c} \includegraphics[width=90mm]{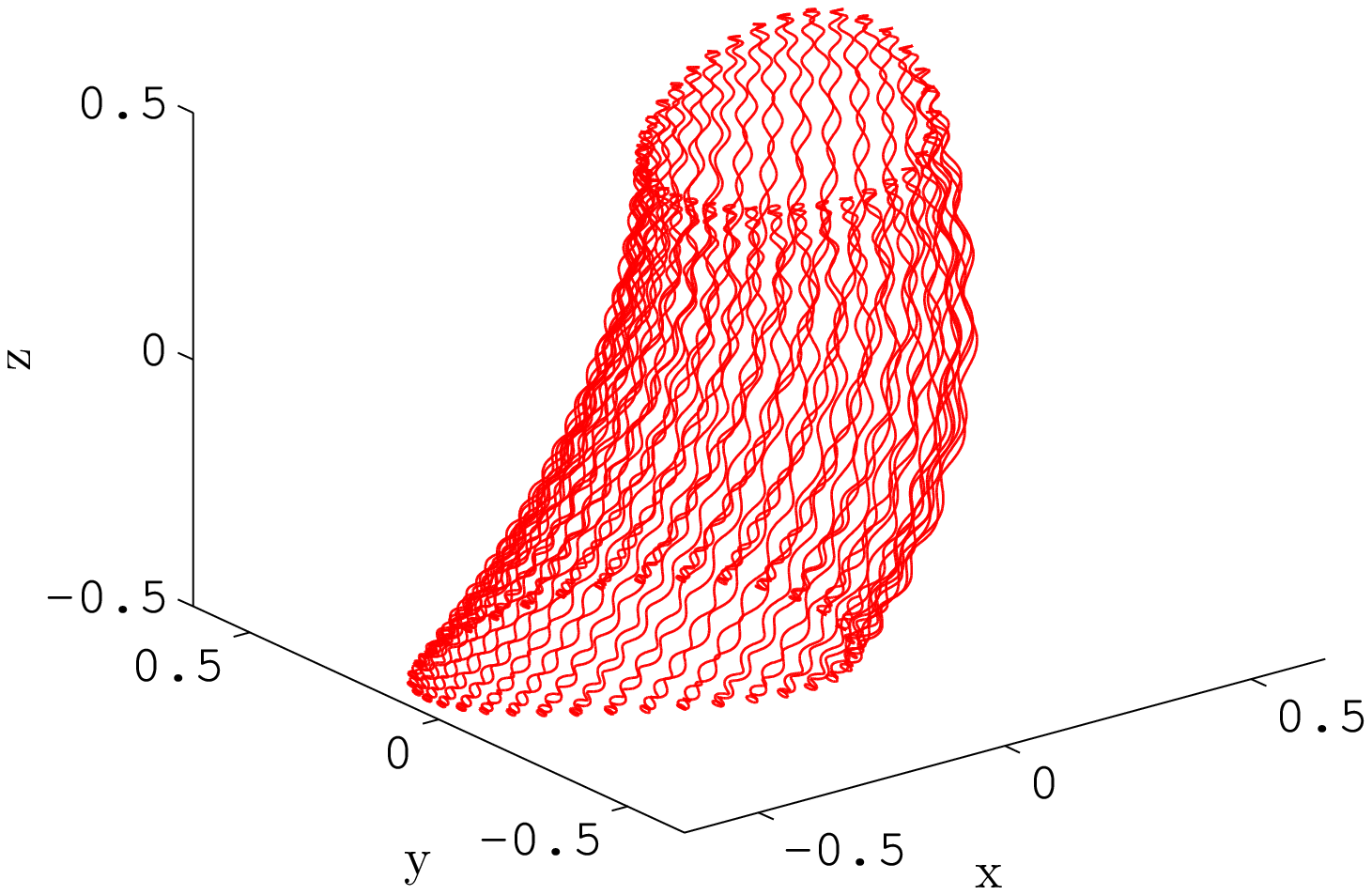}} 
   \caption{Particle trajectories for the three variations of the Penning trap.}
\end{figure}


\subsection{Ideal Penning trap with uniform magnetic field}

The ideal Penning trap consists of a 
static quadrupole electric field and a uniform magnetic field given by
\begin{equation} \label{eq:Penning_fields}
\begin{aligned}
E &= \kappa \begin{bmatrix} x & y & -2z \end{bmatrix}^{\trans}, \\
b &=  \begin{bmatrix} 0 & 0 & b_z \end{bmatrix}^{\trans},
\end{aligned}
\end{equation}
where $\kappa$ depends on the geometry and the voltage of the electrodes. 
The analytic solution for the motion of a single particle in an 
ideal Penning trap can be found in \cite{Kretzschmar}. 
Choosing the trap parameters suitably gives a stable periodic orbit consisting of 
a fast gyromotion perpendicular to $b$ with the frequency $\omega_{mc} \approx \omega_c=\tfrac{c}{m}b_z$, 
a slower axial motion along $b$ with the frequency $\omega_z=\sqrt{2 \frac{c}{m} \kappa}$ solely caused by 
the electrostatic fields and an even slower circular $E \times B$ drift motion around the $z$-axis with the magnetron 
frequency $\omega_m$. 
We set $c=1$, $m=1$, $b_z=100$, $\kappa=10$ and the initial coordinates of the particle  
$q_0=\begin{bmatrix} \tfrac{1}{3} & 0 & \tfrac{1}{2} \end{bmatrix}^{\trans}$ and 
$p_0=\begin{bmatrix} 0 & 1 & 0 \end{bmatrix}^{\trans}$. 
The particle trajectory of this example is shown in Figure~\ref{fig:trajectories_a}.

As the magnetic field is uniform in space the central step of the symmetric Hamiltonian splitting 
\eqref{eq:symm_split} can be computed exactly as is done in Scovel's method. We compare this method
to the method of Chin \eqref{eq:Chinb}, the method of Spreiter and Walter 
(denoted by SpW),
the Boris pusher \eqref{eq:Boris_splitted} and its matrix exponential variation \eqref{eq:Boris_splitted2} 
(denoted by Boris exp).

\begin{figure}
   \centering
      \subfloat[Relative error of energy.]{\label{fig:penning_ideal_energy}
\hspace{-5mm} \includegraphics[width=90mm]{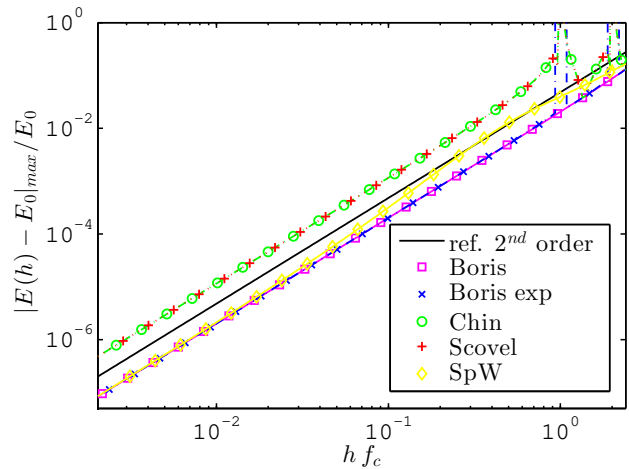}}\\
      \subfloat[Absolute error of position.]{\label{fig:penning_ideal_position}
\hspace{-5mm} \includegraphics[width=90mm]{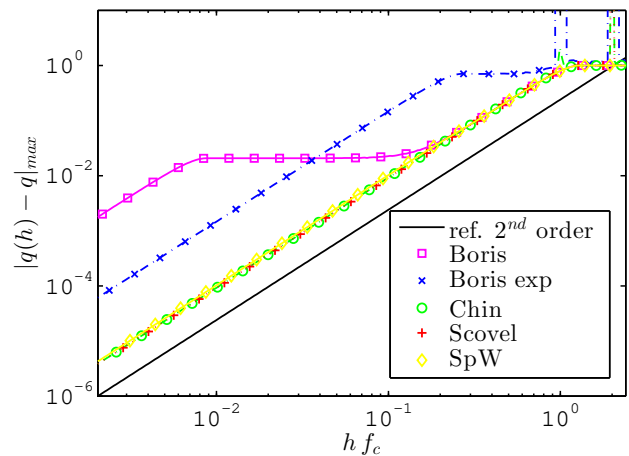}}
   \caption{The ideal Penning trap: maximum error of the total energy $E$ and of the position $q$ 
along a trajectory of two complete magnetron cycles as a function of the time step size $h$ 
(in units of the inverse cyclotron frequency $f_c^{-1}$).}
\end{figure}

We compare the methods by their ability to deliver a stable orbit over a large number 
of magnetron cycles (full circles around the $z$-axis), by the relative energy error 
(Figure~\ref{fig:penning_ideal_energy}) and by the absolute error of the position
(Figure~\ref{fig:penning_ideal_position}). As both types of errors are heavily oscillating,
the maximum value of the error along a trajectory with a total simulation time of
$2 \pi \omega_m^{-1}$ (a complete magnetron cycle) is plotted against the time
step length, using logarithmic scales for both quantities. The time step size ranges from
0.002 to 2.4 cyclotron cycles $f_c^{-1} = \tfrac{2\pi}{\omega_c}$, violating the stability
criteria for some methods. A further reference line with a slope
corresponding to a second-order method has been added to both plots.

With the exception of Spreiter and Walter all methods feature a bounded energy error and a position error
increasing linearly in time on a closed orbit, as long as the time step size is considerably smaller than $f_c^{-1}$.
The method of Spreiter and Walter shows a linearly increasing energy error until its orbit
becomes unstable, so with longer simulation time its errors become worse relative to the other methods.

For larger time steps the method of Chin, the Hamiltonian splitting (Scovel) and the Boris pusher in its matrix
exponential variation become unstable for $f_c h \approx n$, $n$ integer, while the Boris pusher
using the Cayley transform preserves a closed orbit and a bounded energy error.

The Boris pusher provides in both variations  a smaller energy error than the method of Chin and the
Hamiltonian splitting, but shows a larger position error for small time steps.
Figure~\ref{fig:penning_ideal_position} demonstrates an interesting behavior of the Boris pusher:
for a large range of time step sizes the maximum position error remains constant at $\Delta q = 0.02$,
the diameter of the gyromotion. This can be explained by the fact that the Boris pusher can reproduce
the correct first-order gyrocenter drift motion \cite{ParkerBirdsall}.

The position error of a stable trajectory is bounded by 
the maximum diameter of the closed orbit (approximately 1 in our example),
which causes the cut-offs in Figures ~\ref{fig:penning_ideal_position}, ~\ref{fig:penning_bottle_position}, and
~\ref{fig:penning_asymm_position}.

\subsection{Penning trap with magnetic bottle}

The third experiment features a nonuniform magnetic field: a Penning trap with a magnetic bottle 
determined the magnetic field
\begin{equation} \label{eq:magnetic_bottle}
b = b_z \begin{bmatrix} 0 & 0 & 1 \end{bmatrix}^{\trans} + \beta 
\begin{bmatrix} -xz & -yz & z^2-\tfrac{x^2+y^2}{2} \end{bmatrix}^{\trans}.
\end{equation}
Here $E$ is as given in \eqref{eq:Penning_fields}.
A physical implementation of a Penning trap featuring such a magnetic field can be found in \cite{Rodegheri}.
We set $\beta=200$ and all the other parameters are as in the last example. 
The strength of the magnetic field and the cyclotron frequency $\omega_c$ vary only a little along the trajectory. 
The particle trajectory of this example is shown in Figure~\ref{fig:trajectories_b}.

\begin{figure}
   \centering
      \subfloat[Relative error of energy.]{\label{fig:penning_bottle_fp_energy}
\hspace{-5mm}\includegraphics[width=90mm]{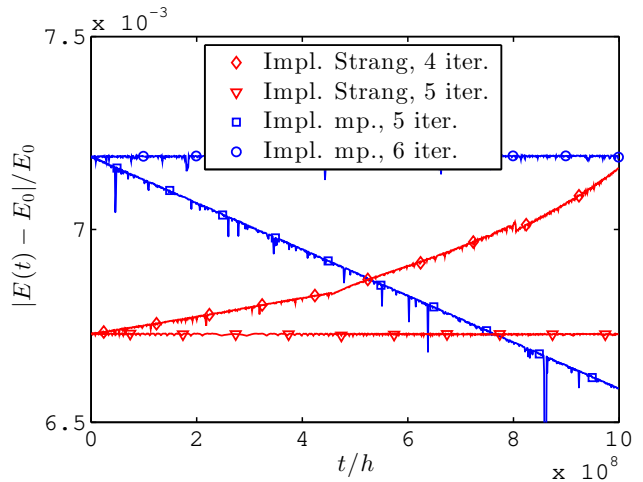}}\\
      \subfloat[Relative error of $\alpha$.]{\label{fig:penning_bottle_fp_phasexy}
\hspace{-5mm}\includegraphics[width=90mm]{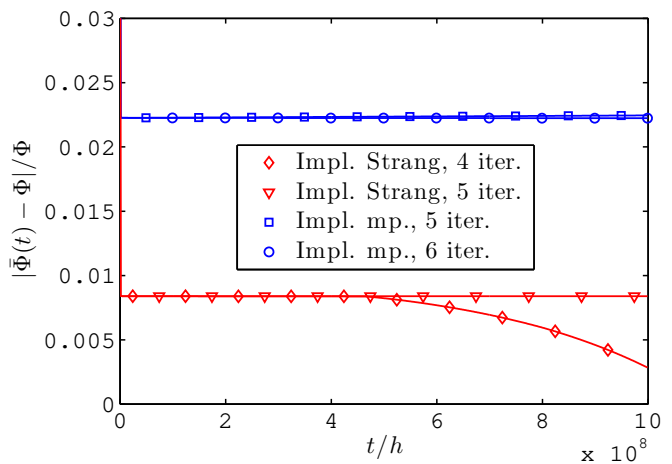}}
   \caption{Penning trap (magnetic bottle): long time evolution of the errors of the 
total energy $E$ and the angle of circular motion $\alpha$ for different numbers of 
fixed-point iterations over $10^9$ time steps of fixed size $h=0.4 f_c^{-1}$.}
\end{figure}

When using the symmetric splitting method \eqref{eq:symm_split}, approximations
have to be made for the mid-step since now the magnetic field is nonuniform in space.
Figure~\ref{fig:penning_bottle_fp_energy} shows the effect of the number of employed fixed-point iterations
on the energy preservation of this method in a long time simulation of a trapped 
particle using a time step size of $h=0.4 \, f_{c,0}^{-1}$, 
where $f_{c,0}$ is the cyclotron frequency at the initial position. 
As the position error is bounded by the closed orbit we choose the relative error of $\alpha$ 
(defined as the
angle of motion around the $z$-axis) as a measure for the trajectory's error in 
Figure~\ref{fig:penning_bottle_fp_phasexy}. Both errors are heavily oscillating, 
so the maximum value over each $10^6$ consecutive time steps is drawn.

The mid-step is either approximated by the implicit method \eqref{eq:subsplit1}
(denoted as Impl.~Strang) or by \eqref{eq:subsplit2} (denoted as Impl.~mp.).
For 5 fixed-point iterations (in the case of Impl.~Strang) and 6 iterations 
(in the case of Impl.~mp.) both alternatives seem to become approximately energy preserving and provide a stable orbit.
This is expected from the theory, since the order of symmetry of these methods equals the number of fixed-point iterations
employed \cite{EinkemmerOstermann} and since a good long time behavior can be expected for 
symmetric methods \cite{HairerLubichWanner}.

We compare the methods of Chin, the Boris pusher and the symmetric splitting method 
\eqref{eq:symm_split} using the implicit schemes
as described above. Furthermore, we apply a high-order composition of \eqref{eq:subsplit1} 
or \eqref{eq:subsplit2} (denoted as Impl.~mp.~comp.).
The composition is performed using an 8th-order symmetric scheme of Suzuki and Umeno \cite{SuzukiUmeno} 
consisting of 15 substeps. Each substep uses 16 fixed-point iterations. Consequently, 
such a composition scheme requires a large number of magnetic field evaluations per time step: 
240 for the method \eqref{eq:subsplit2}, or 255 when using the method \eqref{eq:subsplit1} 
due to the calculation of the initial value within each substep. Nevertheless, 
a single electric field evaluation suffices for each time step.
As the results of the implicit schemes do not differ considerably when the composition 
is applied, only the ones of \eqref{eq:subsplit2}
are drawn in the plots.

As before we provide two log-log plots showing the dependence of the 
energy and the position errors on the time step size
in the range of 0.002 to 0.4 cyclotron cycles $f_{c,0}^{-1}$. Lacking an analytic 
solution to determine the position errors we use the 
trajectory of a high-order method with a small time step (0.00025 cyclotron cycles) 
as a reference: the fourth-order method of \cite{McLachlan} 
applied on $\varphi_h^{\mathcal{T}}$ and $\varphi_h^{\mathcal{E} + \mathcal{B}}$
(denoted as M in \cite{Chin}) in combination with a composition scheme 
\cite{Suzuki_Fractals} increasing the order to 6.

\begin{figure}
   \centering
      \subfloat[Relative error of energy.]{\label{fig:penning_bottle_energy} 
\hspace{-5mm}\includegraphics[width=90mm]{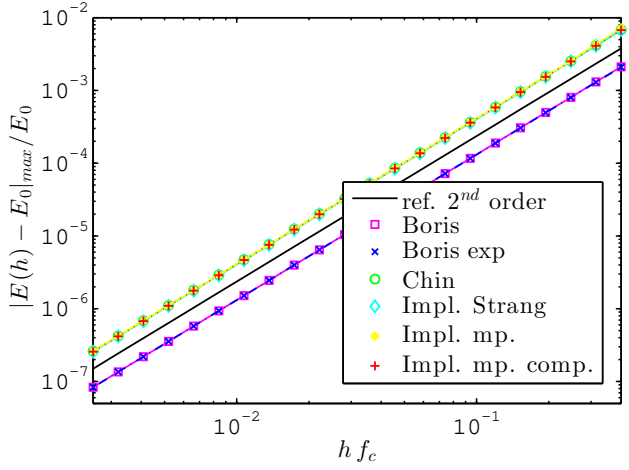}}\\
      \subfloat[Absolute error of position.]{\label{fig:penning_bottle_position} 
\hspace{-5mm}\includegraphics[width=90mm]{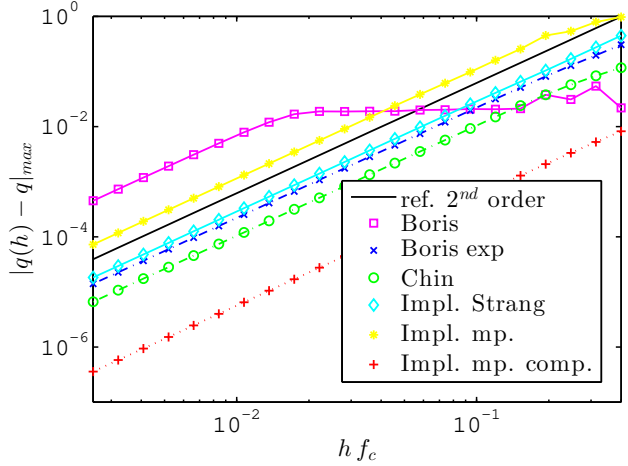}}
   \caption{Penning trap (magnetic bottle): maximum error (along trajectory) of 
the total energy $E$ and of the position $q$ as a function of the time step size $h$ 
(in units of the inverse cyclotron frequency $f_c^{-1}$).}
\end{figure}

All variations of the symmetric Hamiltonian splitting method give an energy error 
similar to Chin's scheme (Figure~\ref{fig:penning_bottle_energy}). 
However, the smallest position error is obtained by using 
the composition scheme for the mid-step of the Hamiltonian splitting.

\subsection{Penning trap with asymmetric magnetic field}

\begin{figure}
   \centering
      \subfloat[Relative error of energy.]{\label{fig:penning_asymm_energy} 
\hspace{-5mm}\includegraphics[width=90mm]{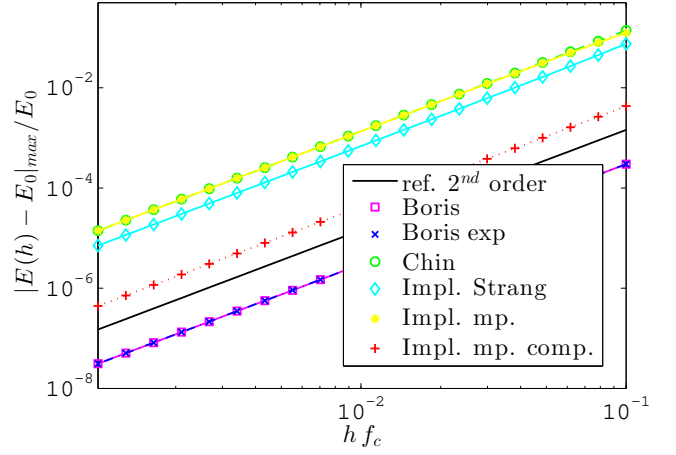}}\\
      \subfloat[Absolute error of position.]{\label{fig:penning_asymm_position} 
\hspace{-5mm}\includegraphics[width=90mm]{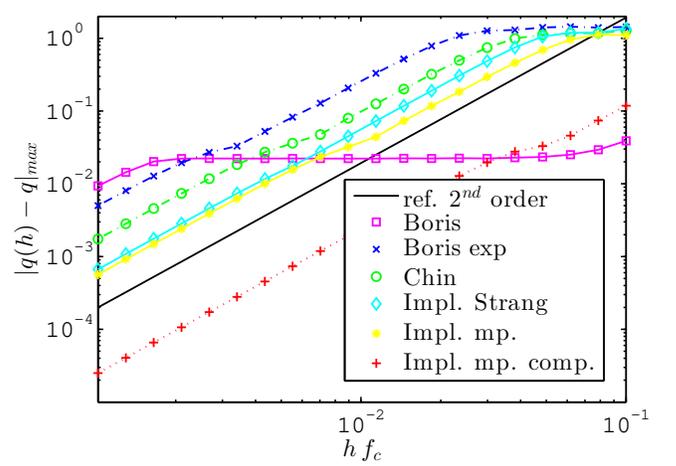}}
   \caption{Penning trap (asymmetric field): maximum error (along trajectory) of the total energy $E$ and of the position $q$ as a function of the time step size $h$ (in units of the inverse cyclotron frequency $f_c^{-1}$).}
\end{figure}

\begin{figure}
   \centering
      \subfloat[Relative error of energy $E$.]{\label{fig:penning_asymm_longrun_energy}
\hspace{-5mm}\includegraphics[width=90mm]{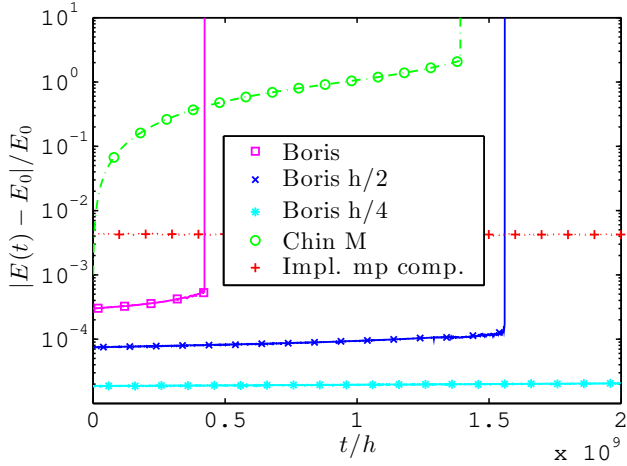}}\\
      \subfloat[Magnetic moment $\mu$.]{\label{fig:penning_asymm_longrun_magmom}
\hspace{-5mm}\includegraphics[width=90mm]{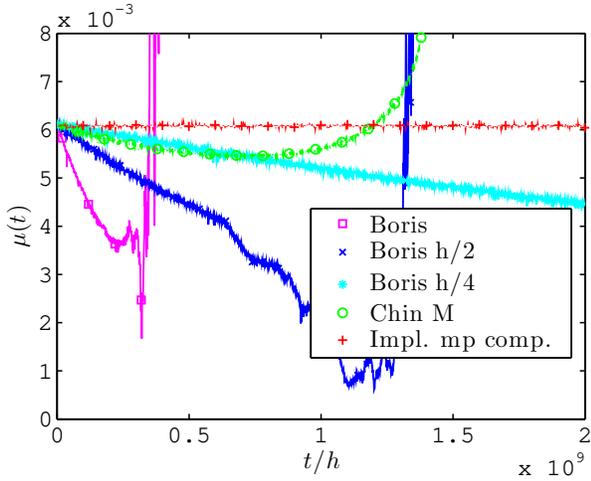}}
   \caption{Penning trap (asymmetric field): long time evolution of the magnetic moment 
$\mu$ and the error of the total energy $E$ over $2 \cdot 10^9$ time steps of fixed size $h=0.1 f_c^{-1}$.}
\end{figure}

As the Penning trap with magnetic bottle is still rotationally invariant, we also consider an example 
with an asymmetric magnetic field given by
\begin{equation} \label{eq:unsym_mag_field}
b(x,y,z) = b_z \begin{bmatrix} \tfrac{1}{3} & 0 & 1 \end{bmatrix}^{\trans} + \beta \begin{bmatrix} y-z & x+z & y-x \end{bmatrix}^{\trans},
\end{equation}
where $\beta=50$. Notice that $\nabla \cdot b = \nabla \times b = 0$.
$E$ is again as given in \eqref{eq:Penning_fields}.
The particle trajectory of this example is shown in Figure~\ref{fig:trajectories_c}.
Along the trajectory the magnetic field strength increases by a factor of 1.8. Therefore time step sizes 
between 0.001 and 0.1$f_{c,0}^{-1}$ are investigated in a simulation lasting for approximately 24 orbit cycles.
Figures~\ref{fig:penning_asymm_energy} and~\ref{fig:penning_asymm_position} show the energy and position errors for the
methods giving similar results as in the previous example. The reference trajectory for the position errors 
has been computed with
the method described in the previous example.
Using the largest time step size $h=0.1 f_{c,0}^{-1}$ shows unbounded energy errors for all the methods with 
the exception of the Boris pusher with the Cayley transform
and the fixed point iteration schemes of the symmetric Hamiltonian splitting in combination with the high-order composition method
of the mid-step. Therefore another long time integration over $2 \cdot 10^9$ time steps of length 0.1$f_{c,0}^{-1}$ 
has been performed using both methods and the fourth-order method of \cite{McLachlan,Chin} (denoted by Chin M).
For this example we measure also the magnetic moment defined by 
$$
\mu(q,p) = \frac{m}{2c^2} \frac{\norm{\Omega(q)p}^2}{ \norm{ b(q) }^3 }  = \frac{c}{2m^2} \frac{\norm{\Omega(q)p}^2}{ \omega_c^3 }.
$$
It is a so-called adiabatic invariant of the system   \cite{HairerLubichWanner} and should be conserved approximately.
The magnetic moment as well as the energy error are oscillating within some range, so the maximum value on consecutive
intervals of $10^6$ time steps is drawn in 
Figures~\ref{fig:penning_asymm_longrun_energy}~and~\ref{fig:penning_asymm_longrun_magmom}.

While the energy error of the Boris pusher remains smaller than that of the symmetric Hamiltonian splitting for
over $4 \cdot 10^8$ time steps, the evolution of the magnetic moment $\mu$ indicates the instability long before. 

An initial decrease of $\mu$ is followed by a rapid growth after $3 \cdot 10^8$ time steps for the Boris scheme. 
The Boris scheme with 2 times smaller time steps (denoted Boris h/2) shows a similar behavior with a delay.
The fourth-order method denoted by Chin M exhibits a quickly growing error in energy but a smoother change of $\mu$.
Only the symmetric Hamiltonian splitting in combination with a high-order approximation of the mid-step is able to provide a stable orbit.

\section{Conclusions}

We have considered here time integrators for particle trajectories in combined electric and magnetic fields.
We have reformulated and analyzed existing methods found from literature under a common framework and additionally
proposed a new method which deals the effect of the magnetic field implicitly.
This is done by splitting the right-hand side of the equations of motion into three terms and by using the matrix function
approach. The matrix functions are used for the magnetic field term, and the so-called
Rodrigues formula simplifies their evaluation. We represent the equations of motion as a Poisson system
and analyze the structure preserving properties of the integrators.

The existing methods we found in literature are: Boris--Buneman scheme, Chin's schemes, 
the method of Spreiter and Walter, and the Scovel's method.
Additionally we have proposed a new implicit scheme. This scheme is the same as 
the Scovel's method when the magnetic field is homogeneous.

The Boris--Buneman scheme and the Chin's schemes are shown to be both symmetric and volume preserving.
For the Spreiter and Walter's scheme we did not find any structure preserving properties.
The new method we propose is fully structure preserving up to an order
equaling the number of fixed-point iterations used for the implicit mid-step.
In numerical experiments we illustrate the effect of the number of fixed-point
iterations on the preservation of the structure. In our experiments 6 iterations were enough
to reduce the error negligible.

As a first numerical comparison, we tested all the schemes except Spreiter and Walter's method 
in a case where there is no electric field. Due to the time symmetry, all the methods seemed to preserve invariants.
For the position error the newly proposed implicit splitting method seemed to be the most accurate
for small time steps, whereas for larger time steps the method of Boris was the most accurate.
When the motion was separated into a periodic and a drift term, it was found that the
method of Boris preserved the drift motion the best and the periodic motion the worst,
whereas this was the opposite for the implicit splitting method.

To compare schemes in combined electric and magnetic fields, we integrated particle trajectories within a Penning trap. 
As the first example of combined fields we considered a constant homogeneous magnetic field, 
so Scovel's method could be used producing energy and solution errors very similar to the method of Chin: 
a lower solution error than the Boris pusher for step-sizes considerably smaller than the inverse cyclotron frequency
at the cost of a slightly larger but likewise bounded energy error. 
Spreiter and Walter's method on the other hand did not feature a bounded energy error, 
so the orbit of the particle was unstable.

As a second Penning trap example we have an inhomogeneous but radially symmetric magnetic field. 
Replacing Scovel's scheme by the new method still gave the same energy error like Chin's scheme. 
The solution error differed depending on the method used for the implicit mid-step iteration. 
Applying a high order composition scheme on the mid-step provided the lowest solution error of all methods tested.

In the last numerical example of a Penning trap, we performed long time integration of a trajectory in an asymmetric
magnetic field and compared three methods: one of the Chin's schemes, the Boris--Buneman scheme
and the newly proposed implicit method. We followed the preservation of the energy and of the magnetic moment.
As expected from the analysis of Section \ref{sec:methods}, the implicit method turned out to be the most stable. 
The Boris--Buneman scheme, being the second most stable, was more unstable than the implicit scheme even for 
4 times smaller time steps.

\section{Appendix: Matrix functions} \label{sec:appendix}

The following identities are easily verified for the skew-symmetric matrix $\Omega$ (see \eqref{def:Omega}):
\begin{equation*} 
\begin{aligned}
\Omega^2        &= \omega \omega^{\trans} - \omega_c^2 I, \\
  \Omega^3        &= - \omega_c^2 \Omega,
\end{aligned}
\end{equation*}
where $\omega_c = \sqrt{\omega_1^2 + \omega_2^2 + \omega_3^2}$.
Using the second identity and the Taylor series of $\sin$ and $\cos$,
we get the Rodrigues formula for the exponential of $h \Omega$:
\begin{equation} \label{eq:Rodrigues}
   \exp{(h\Omega)} = I + \frac{\sin{h\omega_c}}{\omega_c}\Omega + \frac{1-\cos{h\omega_c}}{\omega_c^2}\Omega^2.
\end{equation}

\subsection{$\phi$ functions}

We list here the necessary definitions and properties of the $\phi$ functions.
These functions are defined by
\begin{equation}\label{eq:def-phi}
\phi_0(z) = \ee^z,\qquad \phi_\ell (z) = \int_0^1 \ee^{(1-\theta)z} 
\frac{\theta^{\ell-1}}{(\ell-1)!}\,\dd\theta, \quad \ell\geq 1.
\end{equation}
They satisfy $\phi_\ell(0) = \frac1{\ell!}$ and the recurrence relation
\begin{equation}\label{eq:rec}
\phi_{\ell+1}(z) = \frac{\phi_\ell(z) - \phi_\ell(0)}{z},\quad \ell\ge 0.
\end{equation}
The $\phi$ functions are analytic in the whole complex plane. Therefore,
the matrix $\phi(A)$ is well defined for every $A \in \mathbb{C}^{n \times n}$.
For various definitions and more details on matrix functions we refer to \cite{Higham}.
We mention that the $\phi$ functions play a crucial role also in the so-called exponential integrators \cite{HochbruckOstermann}.

From \eqref{eq:rec} it follows that the series representation 
of the $\phi_\ell$ function is given by
\begin{equation} \label{eq:phi_taylor}
\phi_\ell(z) = \sum_{j=0}^{\infty} \frac{z^j}{(j + \ell)!}.
\end{equation}
We will also need the following formula, which can be easily verified using \eqref{eq:phi_taylor}
\begin{equation} \label{eq:phi_augmented}
  \phi_\ell\left( h \begin{bmatrix} 0 & \frac{1}{m} I \\ 0 & \Omega \end{bmatrix} \right)= 
  \begin{bmatrix} \tfrac{1}{\ell!}I & h \phi_{\ell + 1} ( h \Omega ) \\ 0 & \phi_\ell(h\Omega) \end{bmatrix},\quad \ell\ge 0.
\end{equation}
Using the formula \eqref{eq:Rodrigues} and the integral definition \eqref{eq:def-phi},
the following concise representations can be derived
  \begin{align}
  \phi_1(h\Omega) = I &+ \frac{1 - \cos(h\omega_c)}{h \omega_c^2} \Omega +  
  \frac{h\omega_c - \sin(h\omega_c)}{h \omega_c^3} \Omega^2, \nonumber \\
    \phi_2(h\Omega) = I &+ \frac{h\omega_c - \sin(h\omega_c)}{h^2 \omega_c^3} \Omega \label{eq:Rodrigues_phi} \\
    & + \left( \frac{\cos(h\omega_c) -1}{h^2 \omega_c^4} + \frac{1}{2 \omega_c^2} \right) \Omega^2. \nonumber
  \end{align}

\subsection{Variation-of-constants formula}
The so-called variation-of-constants formula provides the solution of a semilinear
ODE
$$
y'(t) = Ay(t) + g(t,y(t)), \quad y(0) = y_0.
$$
It is given as
\begin{equation} \label{eq:voc}
u(t) = \exp(t A) u_0 + \int_0^t \exp \big( (t-\tau)A \big) \,g(\tau,u(\tau)) \,\dd\tau.
\end{equation}

\subsection{Relations for the Boris pusher}
Let $\Omega$ be as in \eqref{def:Omega}.
An explicit calculation shows that
\begin{equation}  \label{eq:inverse_I-O}
 \big( I - \tfrac{h}{2} \Omega \big)^{-1} = I + \kappa \big( \tfrac{h}{2} \Omega + \tfrac{h^2}{4} \Omega^2 \big),
 \end{equation}
where $\kappa = ( 1+\tfrac{h^2}{4}\omega_c^2 )^{-1}$ and $\omega_c = \sqrt{\omega_1^2 + \omega_2^2 + \omega_3^2} $.
Multiplication by $ I + \tfrac{h}{2} \Omega $ gives 
the Cayley transform of $\tfrac{h}{2} \Omega$,
\begin{equation} \label{rot_Boris} 
\begin{aligned}
  R(\tfrac{h}{2} \Omega) &=
  \big( I - \tfrac{h}{2} \Omega \big)^{-1} \big( I + \tfrac{h}{2} \Omega \big) \\
     &= I + 2 \kappa \big( \tfrac{h}{2} \Omega + \tfrac{h^2}{4} \Omega^2 \big).
\end{aligned}
\end{equation}
The application of $R(\tfrac{h}{2} \Omega)$ on $p^+$ is equivalent to the cross product formulation of the Boris scheme
(see \cite{BirdsallLangdon} for the cross product formulation), since if
 \begin{equation*} 
 \begin{aligned}
 t  &=  \tfrac{h}{2} \omega, \\
 s  &= \tfrac{2t}{1+\norm{t}^2} = \kappa h \omega, \\
 p' &= p^+ + p^+ \times t = \big( I + \tfrac{h}{2} \Omega \big) p^+,
 \end{aligned}
 \end{equation*}
 then
 \begin{equation*}
 \begin{aligned}
 p^{++} &\,\,= p^+ + p' \times s,  \\
 &\,\,=  p^+ + \big( ( I + \tfrac{h}{2} \Omega ) p^+ \big) \times \kappa h \omega, \\
 &\,\,=  p^+ + 2 \kappa \tfrac{h}{2} \Omega \left( I + \tfrac{h}{2} \Omega \right) p^+, \\
 &\,\,=  \big( I + 2 \kappa ( \tfrac{h}{2} \Omega +  \tfrac{h^2}{4} \Omega^2 ) \big) p^+.
  \end{aligned}
 \end{equation*}



\newpage 
%

\end{document}